\documentclass{amsart}[12pt,epsf,amscd,amstex,amssymb]

\newtheorem{theorem}{Theorem}[section]
\newtheorem{corollary}[theorem]{Corollary}
\newtheorem{lemma}[theorem]{Lemma}
\newtheorem{proposition}[theorem]{Proposition}

\newtheorem*{coro}{Corollary} 

\theoremstyle{definition}
\newtheorem{definition}[theorem]{Definition}
\newtheorem{remark}[theorem]{Remark}
\newtheorem{example}[theorem]{Example}

\newcommand{\Ext}{\operatorname{Ext}}
\newcommand{\Hom}{\operatorname{Hom}}
\newcommand{\Ker}{\operatorname{Ker}}
\newcommand{\End}{\operatorname{End}}
\newcommand{\Tor}{\operatorname{Tor}}
\newcommand{\Alc}{\operatorname{Alc}}
\newcommand{\rmod}{\operatorname{mod}}

\newcommand{\tii}{\widetilde}

\newcommand{\kk}{{\bf k}}

\newcommand{\bI}{{\bf I}}
\newcommand{\bLG}{{\bf ^LG}}

\newcommand{\GO}{^L{\bf G}_O}

\newcommand{\bp}{{\bf p}}

\newcommand{\cal}{\mathcal}

\newcommand{\Qlbar}{{\overline {{\mathbb Q}_l} }}

\newcommand{\A}{{\mathcal A}}

\renewcommand{\AA}{{\mathcal A}}

\newcommand{\BB}{{\mathcal B}}
\newcommand{\DD}{{\mathcal D}}
\newcommand{\EE}{{\mathcal E}}
\newcommand{\MM}{{\mathcal M}}
\newcommand{\PP}{{\mathcal P}}
\renewcommand{\SS}{{\mathcal S}}

\newcommand{\OO}{{\mathcal O}}
\newcommand{\E}{{\mathcal E}}

\newcommand{\N}{{\mathcal N}}
\newcommand{\Nt}{{\widetilde{\N}}}
\renewcommand{\gg}{{\mathfrak g}}
\newcommand{\hh}{{\mathfrak h}}
\newcommand{\pp}{{\mathfrak p}}
\newcommand{\gt}{{\tii{\mathfrak g}}}

\newcommand{\de}{\partial}

\newcommand{\tw}{ {}^{(1)} }               

\newcommand{\imbed}{\hookrightarrow}
\newcommand{\iso}{{\tii \longrightarrow}}

\newcommand{\To}{\longrightarrow}

\newcommand{\Gm}{{{\mathbb G}_m}}

\newcommand{\al}{\alpha}

\newcommand{\la}{\lambda}
\newcommand{\La}{\Lambda}
\newcommand{\Ga}{\Gamma}

\newcommand{\oplusl}{\bigoplus\limits}

\newcommand{\bA}{{\bf A}}
\newcommand{\bS}{{\bf S}}
\newcommand{\bY}{{\bf Y}}
\newcommand{\bpr}{{\bf pr}}

\newcommand{\Gr}{{ {\cal G}r}}
\newcommand{\Fl}{{ {\cal F}\ell}}
\newcommand{\Flf}{{\cal{B}}}

\newcommand{\<}{\langle}
\renewcommand{\>}{\rangle}

\newsymbol\ltimes 226E

\renewcommand{\P}{{\mathcal P}}
\newcommand{\Pmixed}{\P^{mix}}

\newcommand{\Ntil}{{\tilde{\mathcal N}}}
\newcommand{\Ntt}{{\tilde{\tilde{\mathcal N}}}}

\newcommand{\Stt}{{\widetilde{St}}}

\newcommand{\h}{{\frak h}}

\renewcommand{\t}{{\frak t}}
\renewcommand{\b}{{\frak b}}

\newcommand{\B}{\Flf}

\newcommand{\gcheck}{{\g\check{\ }}}
\newcommand{\gchecktil}{{\tilde\g\check{\ }}}

\newcommand{\Lg}{{^L{\frak g}}}
\newcommand{\LG}{{^L{G}}}
\newcommand{\LB}{{^L{B}}}

\newcommand{\pit}{{\widetilde \pi}}

\newcommand{\LL}{{\mathcal L}}
\renewcommand{\O}{{\mathcal O}}

\newcommand{\F}{{\mathcal F}}

\newcommand{\ff}{{\mathfrak f}}

\renewcommand{\H}{{\mathcal H}}

\newcommand{\bq}{{\bf q}}

\newcommand{\bbF}{{\mathbb F}}
\newcommand{\bFbar}{{\bar{\Bbb F}}}
\newcommand{\Gcheck}{{G\check{\ }}}

\newcommand{\bG}{{\bf G}}

\newcommand{\Feq}{F}

\newcommand{\Zet}{{\mathbb Z}}
\newcommand{\Pn}{{\mathbb P}^n}

\newcommand{\Atwo}{{\mathbb A}^2}
\newcommand{\Pone}{{\mathbb P}^1}

\newcommand{\Fq}{{\mathbb F}_q}

\newcommand{\bF}{{\mathbb F}}

\newcommand{\Ce}{{\mathbb C}}
\newcommand{\RE}{{\mathbb R}}

\newcommand{\D}{{\mathfrak D}}

\newcommand{\bu}{\bullet}

\renewcommand{\k}{{\bf k}}

\newcommand{\nab}{\nabla}

\newcommand{\del}{\Delta}

\newcommand{\codim}{{\rm co}\dim}

\title[Noncommutative
 Springer Resolution]{Noncommutative
Counterparts of the Springer Resolution}

\author[Roman Bezrukavnikov]{Roman Bezrukavnikov}

\begin{document}

\begin{abstract}
%
Springer resolution of the set of nilpotent elements in a semisimple
Lie algebra plays a central role in geometric representation theory.
A new  structure on this variety has arisen in
several  representation theoretic constructions,
such as the (local) geometric Langlands duality and  modular
representation theory. It is also related to some algebro-geometric
problems, such as the derived equivalence conjecture and description of
 T.~Bridgeland's space of stability
conditions. The structure can be described as a noncommutative
counterpart of the resolution, or as a $t$-structure on the derived
category of the resolution.
 The intriguing 
 fact that the same $t$-structure
appears in these seemingly disparate subjects has strong technical
consequences for  modular representation theory.
\end{abstract}





\maketitle

\tableofcontents

\section{Introduction}

Springer resolution
 of the variety of nilpotent elements in a
semi-simple Lie algebra
 is
ubiquitous in geometric representation theory. In this article we
show that, besides of this well-known resolution of singularities,
the  variety of nilpotents, as well as some other closely related varieties,
admits a particular {\em noncommutative resolution of singularities},
which arises in different representation theoretic and
algebro-geometric constructions. Here by a noncommutative resolution
of a singular variety $Y$ we mean, following, e.g., \cite{vdB}, a
coherent sheaf of associative $\OO_Y$ algebras satisfying certain
natural conditions, and defined up to a Morita equivalence.

 The constructions are related to such subjects as:
the (local) geometric Langlands duality program
and categorification of representation theory of affine Hecke algebras,
  representation theory of modular Lie algebras and quantum enveloping
algebras at roots of unity, Bridgeland's theory of stability conditions
on triangulated categories, and categorical MacKay
correspondence and generalizations.

\medskip

Let  $G$ be a semi-simple adjoint algebraic group, $\gg$ be its Lie
algebra and $\N\subset \gg$ be the variety of nilpotent elements.
Let $\BB$ be the variety of Borel subalgebras in $\gg$, also known
as the flag variety of $G$, and $\Nt=T^*(\BB)$ be the cotangent
bundle to $\BB$. The Springer resolution is the moment map
$\pi:\Nt\to \N$.

 Our noncommutative resolution $A$ of
$\N$ comes with an equivalence
 between the derived category $D(A)$ of modules over $A$ and the derived
 category $D(\Nt)$ of coherent sheaves on  $\Nt$.
Thus $A$ is determined uniquely up to Morita
equivalence
 by the
 $t$-structure on $D(\Nt)$ induced by the equivalence, i.e., by the
 image of the subcategory of $A$ modules in $D(A)$ under the
 equivalence. We will call this $t$-structure the {\em
 exotic $t$-structure} and objects of its heart {\em exotic sheaves}.
 Thus an exotic sheaf is a complex of coherent sheaves on $\Nt$ which corresponds to
 an $A$-module under the equivalence $D(A)\cong D(\Nt)$.

 A closely related  data first appeared in \cite{AB}, which can be
 considered as a contribution to a local version of the geometric Langlands duality
 program \cite{BD},  \cite{GaiICM}, \cite{Edik}.
 A typical result of geometric Langlands duality
   is an equivalence between
 some derived category of constructible sheaves on a variety related
 to $\LG$ bundles on a curve $C$ and derived category of coherent sheaves on a variety
 related to $G$ local systems on $C$; here $G$ and $\LG$ are reductive groups,
  which are dual in the sense of Langlands. In the
 local version of the theory the curve $C$ is a punctured formal
 disc $\mathbb D$. The role of the moduli stack of $\LG$  bundles is played by a
homogeneous space for the group $\LG((t))$, where $\LG((t))$ stands
for the group  of maps from $\mathbb D$ to $\LG$ (also
known as the formal loop 
 group).
  An example of such a homogeneous space is the {\em affine flag variety} $\Fl$ of $\LG$.
   For an appropriate choice of the category of constructible
 sheaves, the variety related to $G$ local systems turns out to be
 $\Nt$, or rather the quotient stack $\Nt/G$ of $\Nt$ by the natural action of $G$.
  An equivalence between the derived category of $G$-equivariant coherent sheaves on $\Nt$
  and a certain triangulated category of constructible sheaves on
  $\Fl$ is proved in \cite{AB}.
The image of the subcategory of perverse sheaves on $\Fl$ under this
equivalence turns out to consist of {\em equivariant exotic
sheaves}, which are closely related to exotic sheaves (see section
\ref{mutasect} below).


Another construction leading to exotic sheaves is related to modular
representation theory.

 In the second half of the
20th century various geometric methods for representation theory of
semi-simple Lie algebras over characteristic zero fields have been
developed. One of the culminating points is the Localization Theorem
\cite{BeBe}, 
\cite{BrylKa}, motivated by a conjecture by Kazhdan and Lusztig,
which provides an equivalence between the category of modules over a
semi-simple Lie algebra $\gg$ with a fixed (integral regular)
central character and the category of $D$-modules on the flag
variety $\BB$.
 In the paper
\cite{BMR}I, motivated by Lusztig's extension \cite{LKthII}II of
Kazhdan-Lusztig conjectures to the modular setting, we provide a
similar result for semi-simple Lie algebras over algebraically
closed fields of positive characteristic. More precisely, we
establish
 a derived localization theorem, which is an equivalence between the derived
category of appropriately defined $D$-modules (called crystalline,
or PD $D$-modules) on a flag variety and the derived category of Lie
algebra  modules, where a part of the center, the so-called
Harish-Chandra center, acts by a fixed character.

Furthermore, in the case of positive characteristic 
there is a close relationship between crystalline $D$-modules on a
smooth variety $X$ and coherent sheaves on the cotangent space
$T^*X$ \cite{BMR}I, \cite{OV}.
The algebra of crystalline differential operators has a huge center
provided by invariant polynomials of the $p$-curvature of a $D$-module.
This allows to view the differential operators
 as a sheaf of algebras on 
 the cotangent bundle.  This algebra turns out to be an Azumaya algebra.
In the case of the flag variety this Azumaya algebra splits on the
formal neighborhood of each Springer fiber. Thus the derived
localization theorem yields a full embedding from the category of
finite dimensional $\gg$-modules with a fixed (integral regular)
action of the Harish-Chandra center into the derived category of
coherent sheaves on $\Nt$. It turns out that the image of this
embedding consists precisely of exotic sheaves with proper support.
A similar relation is expected between exotic sheaves over a field
of characteristic zero and representations of the quantum Kac--De
Concini  enveloping algebra at a root of unity \cite{KDC}, and also
with some class of $\widehat{\Lg}$ modules at the critical level
 (cf. \cite{Kobi} and \cite{Edik_Denis} respectively); here
$\widehat{\Lg}$ stands for the affine Kac-Moody algebra
corresponding to the Langlands dual algebra $\Lg$.

Thus exotic sheaves are related, on the one hand, to perverse sheaves on
the affine flag variety for the dual group, and on the other hand,
to modular Lie algebra representations.
Comparison of these  two connections allows one to
apply the known deep results about weights of Frobenius acting on Ext's
between irreducible perverse sheaves to numerical questions about
modular representations, thereby providing a strategy for a proof of
Lusztig's conjectures from \cite{LKthII}II.
 The conjectures relate the classes of irreducible
$\gg$-modules 
  to elements of the {\em
canonical basis} in the Borel-Moore homology of a Springer fiber;
thus our work provides a categorification of the canonical bases in
(co)homology of Springer fibers. See also Remark \ref{cohotilt}
 for an application
to representations of quantum groups.

I also would like to point out some parallels between exotic sheaves
and objects
 arising in the work of algebraic geometers studying
derived categories of coherent sheaves on algebraic varieties.
Exotic sheaves can be described in terms of a certain action of the
affine braid group $B_{aff}$ of $\LG$ on $D(\Nt)$. This description
can be reformulated in terms
 a  map from the set of alcoves (connected components of
the complement to affine coroot hyperplanes in the dual space to the
Cartan algebra of $\gg$ over $\RE$) to the set of $t$-structures on
$D(\Nt)$. 
 A similar data has been used by
Bridgeland in \cite{BrK3} to construct a component in
the space of stability conditions \cite{Brstab},
on the derived categories of coherent sheaves on certain varieties.
See also Examples \ref{flop}, \ref{subr} below.

The appearance of the affine braid group, which can be interpreted
as the fundamental group of the set of regular semisimple conjugacy
classes in the dual group $\LG(\Ce)$, suggests a possibility that
the structures described above admit a natural interpretation via
homological mirror duality, which would identify our derived
category of coherent sheaves with a certain Fukaya type category,
where the action of the affine braid group arises from monodromy of
some family over the space of  regular semisimple conjugacy classes
in $\LG(\Ce)$.

Another connection to algebraic geometry is provided by \cite{BeKal}
and \cite{Kal}. As has been noted above, the derived localization
theorem can be interpreted as a construction of a noncommutative
resolution of the nilpotent cone $\N$ using crystalline differential
operators in positive characteristic. It turns out that for more
general resolutions of singularities, which carry an algebraic
symplectic form, a non-commutative resolution can be constructed by
a similar procedure. The construction
 involves
quantizing the algebraic symplectic variety in characteristic
$p$, and relating modules
over the quantization to coherent sheaves.
It has been carried out in
\cite{BeKal}
 for crepant resolutions of  quotients
$V/\Gamma$, where $V$ is a vector space equipped with a symplectic
form, and $\Gamma$ is a finite subgroup in $Sp(V)$; this yields a
particular case of the so-called {\em categorical MacKay
correspondence}. The particular case when $\Gamma$ is the symmetric
group on $n$ letters acting on $V=(\Atwo)^n$ is related to representations
of the rational Cherednik algebra \cite{BFG}.
In Kaledin's work \cite{Kal}
 the construction is generalized to more general
symplectic resolutions of singularities.

\medskip

In the remainder of the text  we explain some of these contexts  (in
the order which is roughly inverse to the above) in some detail.

This is text is a mixture of an exposition of published results and
announcement of yet unpublished ones; statements for which no reference
is provided, and which are not well-known, are to appear in a future
publication.


{\bf{Notations and conventions}} Throughout the text we work over an
algebraically closed field $\kk$; when a semi-simple group $G$ is
involved, we assume that characteristic of $\kk$ is zero or exceeds
the Coxeter number of $G$.

For an algebraic variety $X$ we let $\O_X$ denote the structure
sheaf, and $D(X)=D^b(Coh_X)$ be the bounded derived category of
coherent sheaves on $X$. Given an action of an algebraic group $H$
on $X$ we write $Coh^H(X)$ for the category of $H$-equivariant
coherent sheaves; given a  coherent
sheaf of associative $\O_X$ algebras  we let $Coh(X,\A)$ be the
 category of sheaves of coherent $\A$ modules; if $\A$ is $H$-equivariant
for an algebraic group $H$ acting on $X$, we let $Coh^H(X,\A)$ be
the category of $H$-equivariant sheaves of coherent $\A$-modules.
We write $D(X)$, $D^H(X)$, $D(\A)$, $D^H(\A)$ for the bounded
derived category of $Coh(X)$, $Coh^H(X)$, $Coh(X,\A)$, $Coh^H(X,\A)$
respectively, and $K(X)$, $K^H(X)$, $K(\A)$, $K^H(\A)$ for the
corresponding Grothendieck groups. In particular, there notations
apply for an algebra $A$ finite over the center of finite type.

%
 The
functors of pull-back, push-forward etc. between categories of
sheaves are understood to be the derived functors.

\medskip

{\bf Acknowledgements} I thank my coauthors S.~Arkhipov,
V.~Ginzburg, D.~Kaledin, I.~Mirkovi\' c, V.~Ostrik and D.~Rumynin
for their contribution to the joint results, without which they
would have never been accomplished.
%
The project described here was conceived during IAS Special Year in
Representation Theory (98/99) led by G.~Lusztig. Most of the results
have been obtained by unraveling the formulas in Lusztig's papers,
thus they owe their existence to him. I have learned Lusztig's
results and many other things from M.~Finkelberg. I have also
benefitted a lot from ideas of I.~Mirkovi\' c and his generosity in
sharing them. Conversation with many people were very helpful, the
incomplete list includes  A.~Beilinson, V.~Drinfeld, D.~Gaitsgory,
V.~Ginzburg, V.~Ostrik.
I am very grateful to all these people.
Finally, I thank M.~Finkelberg, D.~Kazhdan and I.~Mirkovi\' c for
reading a preliminary version of this text and making helpful
comments and suggestions.

\section{Noncommutative resolutions and braid group actions}
\label{se2}

\subsection{Braid group actions and noncommutative Springer resolution}

Though the motivation for the study of our main object comes from
applications to representation theory, we first describe it in the
language of  algebraic geometry. We briefly recall some ideas
 from \cite{BonOr}, \cite{Bridge}, \cite{vdB}.

Let $Z$ be a singular algebraic
 variety. We refer, e.g., to
\cite{vdB} for the notion of a {\em crepant} resolution;
 it is easy to see that resolutions $\pi$, $\pit$ described above
 are crepant.

  By a {\em noncommutative resolution}
\cite{vdB}\footnote{The definition in {\em loc. cit.} is wider, we
use a version convenient for our exposition.}
one means a coherent torsion free sheaf $A$ of associative $\O_Z$
algebras, which is generically a sheaf of matrix algebras and
has finite homological dimension. There exists also a notion
of a noncommutative crepant resolution, see \cite{vdB}.
It has been conjectured in {\em loc. cit.} that any two crepant
 resolutions, commutative or not, are derived equivalent, in particular,
for any crepant resolution $X\to Z$ and any noncommutative
resolution $A$ of $Y$ we have an
equivalence $D(X)\cong D(A)$. 


\subsubsection{The set-up}
Notations $G$, $\gg$, $\BB$, $\pi:\Nt\to \N$ has been defined in the
Introduction. Recall that $\Nt=T^*(\BB)$ parametrizes pairs
$(\b,x)$, where $\b\in \BB$ is a Borel subalgebra, and $x$ is the
element in the nilpotent radical of $\b$. The Springer map
$\pi:\Nt\to \N$ is given by $\pi:(\b,x)\mapsto x$.
 It is embedded in the
{\em Grothendieck simultaneous resolution} $\pit: \gt\to \gg$, where
$\gt$ is the variety of pairs $(\b,x)$, $\b \in \BB$, $x\in \b$, and
$\pit:(\b,x)\mapsto x$. The variety $\gt$ is smooth, and the map
$\pit$ is proper and generically finite of degree $|W|$, where $W$
is the Weyl group. It factors as a composition of a resolution of
singularities $\pit': \gt\to \gg\times_{\hh/W} \hh$ and the finite
projection $\gg\times_{\hh/W} \hh\to \gg$; here $\hh$ is the Cartan
algebra of $\gg$.
 Let $\gg^{reg}\subset \gg$ denote the subspace of regular
(not necessarily semi-simple) elements, and $\gt^{reg}$ be the
preimage of $\gg^{reg}$ in $\gt$; then $\pit'$ induces an
isomorphism $\gt^{reg} \cong \gg^{reg}\times_{\hh/W}\hh$.

 Much of representation theory of $G$ or $\gg$ is in one way or
another related to the geometry of these spaces and maps.

\subsubsection{Affine braid group action}
For a characterization of our noncommutative resolution
we need to introduce some notations.

Let $\Lambda$ be the root lattice of $G$.
 For $\la\in \La$ we will write $\O(\la)$ for the
corresponding $G$-equivariant line  bundle on $\BB$, and we set
$\F(\la)=\F\otimes _{\OO_{\BB}}\O(\la)$ if $\F\in D(X)$ for some $X$
mapping to $\BB$.

Let $W$ be the Weyl group, and set $W_{aff}=W\ltimes \La$.
Then $W$, $W_{aff}$
 are Coxeter groups. Notice that
$W_{aff}$ is the affine Weyl group of the
{\em Langlands dual} group $^LG$.

It was mentioned above that $\gt^{reg}\cong \gg\times_{\hh/W}\hh$;
thus $W$ acts on this space via its action on the second factor.
The formulas $\La\owns \la:\F \mapsto \F(\la)$, $W\owns w:\F \mapsto
w^*(\F)$ are easily shown to define an action\footnote{Throughout
the paper by an action of a group on a category I mean  {\em a weak
action}, i.e., a homomorphism to the group of isomorphism classes of
autoequivalences. I believe that in all the examples in this text a
finer structure can be established, though I have not studied this
question.} of $W_{aff}$ on the category of coherent sheaves on
$\gt^{reg}$.

The  characterization of our ``noncommutative Springer resolution''
  relies on the possibility to extend
this action to a weaker structure on the whole of  $\gt$. To
describe this weaker structure recall that to each Coxeter group
one can associate an Artin braid group; let
$B_{aff}$ denote the group corresponding to $W_{aff}$.
 It admits a  topological interpretation, as the fundamental group
of the space of regular semi-simple conjugacy classes in the
universal cover of the dual group $\LG(\Ce)$.
For $w\in W_{aff}$ consider the minimal decomposition of $w$ as a product
of simple reflection, and take the product of corresponding generators
of $B_{aff}$. This product is well known to be independent on the choice
of the decomposition of $w$, thus we get a map $W_{aff} \to B_{aff}$
which is one-sided inverse  to the canonical surjection $B_{aff}
\to W_{aff}$. We denote this map by $w\mapsto \tii w$. The map is not a
homomorphism, however, we have  $\tii{uv}=\tii u \cdot \tii
v$ for any $u,v\in W_{aff}$ such that $\ell(uv)=\ell(u)+\ell(v)$,
where $\ell(w)$ denotes the length of the minimal decomposition of
$w$.
%
 Let $B_{aff}^+\subset B_{aff}$ be the sub-monoid
generated by $\tii w$, $w\in W_{aff}$.

For a simple reflection $s_\al\in W$ let $\bS
_\al\subset \gt^2$ be the closure
of the graph of $s_\al$ acting on $\gt^{reg}$. We let $S_\al$
denote the intersection of $\bS_\al$ with $\Nt^2$.
Let $\bpr^\al_i:\bS\to \gt$, $pr^\al_i:S\to \Nt$, where $i=1,2$,
be the projections.

Let $\La^+\subset \La$ be the set of dominant weights in $\La$.

For a scheme $Y$ over $\gg$ we set $\tii Y=\Nt\times _\gg Y$, $\tii
\bY=\gt \times _\gg Y$.

\begin{theorem}\label{bract}
a) There exists an (obviously unique) action of $B_{aff}'$
on $D(\gt)$, $D(\Nt)$ such that
 for $\la \in \La^+\subset \La \subset W_{aff}'$
we have $\tii \la:\F \mapsto \F(\la)$
and for a simple reflection
$ s_\al\in W$ we have $\tii s_\al:\F \mapsto (\bpr_1^\al)^*(\bpr_2^\al)_*
\F$ (respectively,
 $\tii s_\al:\F \mapsto (pr_1^\al)^*(pr_2^\al)_*
\F$).

b) This action induces an action on $D(\tii \bY)$, $D(\tii Y)$ for
any scheme $Y$ over $\gg$ such that $\Tor_i^{\O_\gg}(\O_\gt,
\O_Y)=0$, respectively $\Tor_i^{\O_\gg}(\O_\Nt, \O_Y)=0$, for $i>0$.
\end{theorem}

{\em Comment on the proof.}
The Theorem can be deduced from material of either section
 \ref{se3} or \ref{se5} below. The most direct proof relies on
the result of \cite{Rina}.

\begin{remark}\label{tr_slice}
 An example of $Y$ satisfying the assumptions of the
Theorem is given by a transversal slice to a nilpotent orbit. In
particular, if $Y$ is a transversal slice to a subregular orbit,
then $\Nt\times _\gg Y$ is well known to be the minimal resolution
of a simple surface singularity. The affine braid group action in
this case coincides with the one constructed by Bridgeland in
\cite{BrK3}III.
\end{remark}

\begin{remark}\label{heck}
The induced action of $B_{aff}$ on the Grothendieck group
$K(\Nt)$ factors through $W_{aff}$. If one passes to the category
of sheaves equivariant with respect to the multiplicative group,
acting by dilations in the fibers of the projection
$\Nt\to \BB$, then the induced action factors through the {\em
affine Hecke algebra} $\H$, cf. discussion after Theorem \ref{equiv}.
Furthermore, this construction
yields an action of $\H$ on the Grothendieck group
$K(\pi^{-1}(e))$ for each $e\in \N$; these $\H$ modules are called
the standard $\H$-modules.
Thus the Theorem provides a {\em categorification} of the
standard modules for the affine Hecke algebra.
\end{remark}

The next result, which plays an important technical role in the
proofs, is a categorical counterpart of the quadratic relation in
the affine Hecke algebra, see discussion after Theorem \ref{equiv}.

\begin{proposition}
For every simple reflection $s_\al\in W_{aff}$
and every $\F\in \D$ we have a (canonical) isomorphism in the quotient
category $\D/\langle \F\rangle$
$$\tii s_\alpha(\F)\cong \tii s_\alpha^{-1}(\F)\mod \langle \F\rangle.$$
Here $\langle \F\rangle$ denotes the full triangulated subcategory
generated by $\F$.
\end{proposition}


\subsubsection{The $t$-structure and the noncommutative resolution}

We will describe certain noncommutative resolutions $A$, $\bA$ of
$\N$, $\gg\times_{\hh/W} \hh$ respectively, together with
equivalences $D(A)\cong D(\Nt)$, $D(\bA)\cong D(\gt)$, and show how
they appear in representation theory.  Such data is uniquely
determined by the $t$-structures on $D(\Nt)$, $D(\gt)$, which are
the images of the tautological $t$-structures on $D(A)$, $D(\bA)$.

\begin{definition} Let $D$ be a triangulated category equipped with an
action of $B_{aff}$.  A $t$-structure $(D^{<0}, \, D^{\geq 0})$ on
$D$ will be called {\em braid right exact} if any
$b \in B_{aff}^+$ sends $D^{<0}$ to $D^{<0}$.
\end{definition}

\begin{theorem}\label{tstr}
a) Let $X$ be either $\tii Y$ or $\tii \bY$, where $Y\to \gg$ is as
in  Theorem \ref{bract}.

The category $D(X)$ admits a  unique
$t$-structure which is

i) braid right exact, and

ii)  compatible with the standard $t$-structure on
the derived category of vector spaces under the functor of derived global
sections $R\Gamma$.


b) There exists a vector bundle $\E_X$ on $X$, such that the functor
$\F\mapsto R\Hom(\E, \F)$ is an equivalence between $D(X)$ and
$D(A_X)$, sending the $t$-structure described in (a) to the
tautological $t$-structure on $D(A_X)$;  here $A_X=\End(\E_X)^{op}$,
where the upper index denotes the opposite ring.

 Moreover, there exists a
vector bundle $\E=\E_\gt$ on $\gt$, such that for any $X$
we can take $\E_X$ to be the
pull-back of $\E$ to $X$.

\end{theorem}


\begin{remark}
It is clear from the definitions that if $X$ is smooth, then $A_X$
is a noncommutative resolution of $Y\times _\gg \N$ or $Y\times
_{\hh/W}\hh$. In particular, for $Y=\gg$ we get
 $A=A_\Nt$, $\bA = A_{\gt}$, which
are the promised noncommutative resolutions of $\N$,
$\gg\times_{\hh/W}\hh$.
\end{remark}

We will call the $t$-structures described in Theorem \ref{tstr} the
{\em exotic} $t$-structures, the objects of their heart will be
called exotic sheaves.


\begin{example}\label{flop}
Let $G=SL(2)$, thus $\Nt$ is the total space of the line bundle
$\O(-2)$ on $\Pone$, and $\gt$ is the total space of the vector
bundle $\O_\Pone (-1)\oplus \O_\Pone (-1)$. In this case we can set
$\E\cong \O_\gt\oplus \O_\gt(1)$.

This $t$-structure on $D(\gt)$ appeared in Bridgeland's
 proof of the derived equivalence conjecture for varieties of
dimension three  \cite{Brflop}. More precisely, for a flop of
three-folds $X,X'\mapsto Y$ Bridgeland constructs some
noncommutative resolution
 of $Y$ which is derived equivalent to
both $X$ and $X'$. The simplest example of a three-fold flop is as
follows: $X=X'=\gt$, $Y=\gg\times _{\hh/W}\hh$ and the two maps
$X,\, X'\to Y$ are $\pit'$ and $\pit''=\iota \circ \pit'$, where
$\iota$ is an involution of $\gg\times_{\hh/W} \hh$ given by
$(x,h)\mapsto (x,-h)$. The $t$-structure on $D(\gt)$ given by
Bridgeland's construction applied to this flop turns out to coincide
with the $t$-structure provided by Theorem \ref{tstr}.

\end{example}

\begin{example}\label{subr}
Let $Y$ be a transversal slice to the subregular orbit.  Thus $Y$ is
isomorphic to the quotient $\Atwo/\Gamma$ for some finite subgroup
$\Gamma\subset SL(2)$. The fiber product $X=\Nt\times_\gg Y$ is the
minimal resolution of $Y$. It is well known that there exists a
natural equivalence $D(X)\cong D^\Gamma(\Atwo)$. The exotic
$t$-structure coincides with the one induced from the tautological
$t$-structure on $D^\Gamma (\Atwo)$. Thus $A_X$ is Morita equivalent
to the smash product algebra $\Gamma \# \O(\Atwo)$. This
$t$-structure appears also in \cite{BrK3}II.
\end{example}





\subsubsection{Parabolic version}\label{parvar}
One can also consider the partial flag varieties $\PP=G/P$, where
$P\subset G$ is a parabolic subgroup; thus $\PP$ parameterizes
parabolic subalgebras $\pp\subset \gg$ of a given type. There exist
parabolic versions of the Grothendieck-Springer spaces:
$\gt_\PP=\{\pp\in \PP, \, x\in \pp \}$ and $\Nt_\PP=T^*(\PP)$. We
have a proper map $\pi_\PP:\gt \to \gt_\PP$, $(gB,x)\mapsto (gP,x)$.
Also, the projection $G/B\to G/P$
induces a closed embedding $\iota_\PP
:\BB\times _\PP \Nt_\PP\imbed \Nt$; we let $pr_\BB^\PP$ denote the
projection $\BB\times _\PP \Nt_\PP \to \Nt_\PP$.

The following result easily follows from the results of
\cite{BMR}II.

\begin{theorem}
a) There exists a unique  $t$-structure on $D(\gt_\PP)$, whose heart
contains the image of exotic sheaves under the functor $R\pi_{\PP*}:
D(\gt) \to D(\gt_\PP)$. 

b) There exists a unique  $t$-structure on $D(\Nt_\PP)$, such that
for any object $\F$ in its heart the object
$(\iota_{\PP*}(pr^\BB_\PP)^* \F)(\rho)$ is an exotic sheaf.
\end{theorem}

One also has induced nice $t$-structures on $D(Y\times _\gg
\gt_\PP)$, $D(Y\times _\gg \Nt_\PP)$ for $Y$ satisfying a Tor
vanishing condition; we omit the details to save space.

\begin{example}
Let $G=SL(n+1)$ and $\PP=\Pn$. The heart of the $t$-structure on
$\Nt_\PP=T^*\Pn$ has a projective generator $\oplusl_{i=0}^n
\O_{T^*\Pn}(-i)$. The heart of the $t$-structure on $\gt_{\Pn}$ has
a projective generator $\oplusl_{i=0}^n \O_{\gt_\PP}(i)$.
 \end{example}

\subsubsection{Reformulation in terms of $t$-structure assigned to
alcoves} 
A connected component of the complement to the coroot
hyperplanes $H_\al$ 
 in the dual space to {\em real} Cartan algebra
$\hh^*_\RE$ is called an alcove; in particular, the {\em fundamental
alcove} $A_0$ is the locus of points where all positive coroots take
value between zero and one. Let $\Alc$ be the set of alcoves. For
$A_1, \, A_2\in \Alc$ we will say that $A_1$ lies above $A_2$ if for
any positive coroot $\check\alpha$ and $n\in \Zet$, such that the
affine hyperplane $H_{\check\al,n}=\{\la,\ |\, \langle \check
\alpha, \la\rangle =n \}$
 separates $A_1$ and $A_2$, $A_1$ lies above $H_{\al,n}$, while $A_2$ lies
below $H_{\al,n}$, i.e. for $\mu\in A_2$, $\la\in A_1$ we have $\<
\al,\mu\> < n < \<\al,\la\>$ .

\begin{lemma}
There exists a unique map $\Alc\times \Alc\to B_{aff}$,
 $(A_1,A_2)\mapsto b_{A_1,A_2}$, such that

i) $b_{A_2, A_3} b_{A_1,A_2}=b_{A_1A_3}$ for any $A_1,A_2,A_3\in \Alc$.

ii) $b_{A_1,A_2}=\tii w$, provided that $A_2$
lies above $A_1$. Here $w\in W_{aff}$ is such that $w(A_1)=A_2$.
\end{lemma}



The following result is equivalent to Theorem  \ref{tstr}.


\begin{theorem}\label{tstralc}
Let $X=\tii Y$ or $\tii \bY$, where $Y$ is as in Theorem
\ref{bract}.

There exists a unique collection of $t$-structures indexed by
alcoves, $(D^{\leq 0}_A(X),\, D^{>0}_A(X))$  such that:

0) (Normalization)
 The derived global sections  functor $R\Gamma$ is $t$-exact with
respect to the $t$-structure corresponding to $A_0$.

1) (Compatibility with the braid action)
 The action of the element $b_{A_1,A_2}$ sends the $t$-structure
corresponding to $A_1$ to the $t$-structure corresponding to $A_2$.

2) (Monotonicity)
 If $A_1$ lies above $A_2$, then $D^{>0}_{A_1}(X)\supset D^{>0}_{A_2}(X)$.


\end{theorem}

\begin{remark}
The exotic $t$-structure described in Theorem \ref{tstr}
is the one attached to the fundamental alcove $A_0$ by the construction of
 Theorem \ref{tstralc}.
\end{remark}

\begin{remark}
The data described in Theorem \ref{tstralc} resembles the one
obtained by Bridgeland in the course of description of the manifold
of stability conditions on some derived categories of coherent
sheaves. To enhance this point we mention a positivity property of
the $t$-structure $ (D^{\leq 0}_A(X),\, D^{>0}_A(X))$; such
properties play a role in the definition of stability conditions
\cite{Brstab}.
\end{remark}

It is easy to show that 
each of the above $t$-structures  induces a $t$-structure on the
full subcategory $D^f(X)\subset D(X)$ consisting of complexes whose
cohomology sheaves have proper support. Let $\A_A = D^{\leq 0}_A(X)
\cap D^{\geq 0}_A(X)$ be the heart of the $t$-structure, and set
$\A_A^f=\A_A\cap D^f(X)$. It is easy to show that $\A_A^f$ consists
of objects {\em of finite length} in $\A_A$.

Assume that $\kk=\Ce$ and $X$ is smooth. Recall that for a smooth
complex variety $X$
 we have the Chern character map $K(D^f(X))\to
H_*^{BM}(X)$, where $H_*^{BM}$ stands for the Borel-Moore homology
of the corresponding complex variety endowed with the classical
topology.  We have a perfect pairing between cohomology and
Borel-Moore homology.

We have a well-known identification $\hh^*
=H^2(\BB)$.
\begin{proposition}\label{pospos}
For $A\in \Alc$, $\F\in \A_A^f$, $\F\ne 0$ and $x\in A\subset
\hh^*_\RE \subset  H^2(\BB)$ we have
$$\< ch(\F), pr^* (\exp (x)) \> >0,$$
where $exp(x)=1+x+\frac{x^2}{2}+\cdots +\frac{x^{\dim \BB}}{(\dim
\BB)!}$, and $pr$ stands for the projection $X\to \BB$.
\end{proposition}

Finally, we describe compatibility of our $t$-structures with duality.

Let $\SS$ denote the Grothendieck-Serre duality functor.

\begin{proposition}\label{minusA}
$\SS$ sends $\A_A^f$ to $\A_{-A}^f$, where $-A$ denotes the alcove
opposite to $A$.
\end{proposition}

\subsection{Equivariant category and mutations of exceptional
sets}\label{mutasect}
 The categories $D(\Nt)$, $D(\gt)$ have
equivariant versions $D^G(\Nt)$, $D^G(\gt)$. It turns out that these
equivariant categories carry  $t$-structures which are, on the one
hand, closely related to the above $t$-structures on non-equivariant
categories, and, on the other hand, admit a direct description in
terms of generating exceptional sets in a triangulated category.

{\em Until the end of \ref{neqgr} we assume that $char(\kk)=0$.}

\subsubsection{Exceptional sets and mutations}
Recall that an ordered set of objects  $\nab=\{ \nab^i,\;  i\in I\}$
in a triangulated category is called {\it exceptional} if we have
$\Hom^\bu (\nab^i, \nab^j)=0$ for $i<j$;
 $\Hom ^{n}(\nab^i,\nab^i)=0$ for $n\ne 0$,
 and $\End(\nab^i)=\k$.
A set  $\del=\{ \del_i,\;  i\in I\}$ of objects is called {\em dual}
to  $\nab$ if $\Hom^\bu(\del_i,\nab^i)=\k$, and
$\Hom^\bu(\del_i,\nab^j)=0$ for $i\ne j$; it is exceptional provided
$\nab$ is, where the order on $\del$ is defined to be opposite to
that on $\nab$.
 Let  $\nab$, $\del$ be two dual
 exceptional sets which generate a triangulated category $\D$; assume that $\{j\ |\
j\leq i\}$ is finite for every $i\in I$. Then  there exists a unique
$t$-structure $(\D^{\geq 0}, \D^{<0})$ on $\D$, such that
$\nab\subset \D^{\geq 0}$; $\del\subset \D^{\leq 0}$. This
construction is closely related to the definition of a perverse
sheaf, see \cite{Hum} for details.

%
%
 Let $(I,\preceq)$
be an ordered set, and $\nab^i\in \D$, $i\in I$ be an exceptional
set. Let $\leq$ be another order on $I$; we assume that $\{j\ |\
j\leq i\}$ is finite for every $i\in I$. We let $\D_{\leq i}$ be the
full triangulated subcategory generated by  $ \nab^j$, $j\leq i$,
and similarly for  $\D_{< i}$. Then for $i\in I$ there exists a
unique (up to a unique isomorphism) object $\nab^i_{mut}$ such that
$\nab^i_{mut}\in \D_{\leq i}\cap \D_{<i}^\perp$, and
$\nab^i_{mut}\cong \nab^i\mod \D_{<i}$ (see e.g. \cite{Hum}).
 The objects $\nab^i_{mut}$ form an exceptional set indexed by
$(I,\leq)$.


 We will say that the exceptional set $(\nab^i_{mut})$ is
 the $ \leq$ mutation of $(\nab^i)$. This construction is
related, cf. \cite{Hum}, to the action of the braid group on the set
of exceptional sets in a given triangulated category constructed in
\cite{BoKa}, this action is also called the action  by mutations.

\subsubsection{Exceptional sets in $D^G(\Nt)$}
Recall the standard partial order $\preceq$ on the set $\La$ of
weights of $G$, which is given by: $\la \preceq \mu$ if $\mu-\la$ is
a sum of positive roots. Then
 line bundles $\O_\Nt(\la)$ generate $D^G(\Nt)$, and
we have $\Hom^\bu(\O(\la), \O(\mu))=0$ unless $\mu \preceq \la$ and
$\Hom^\bu(\O(\la), \O(\la))=\k$ \cite{Hum}. Thus   for any complete
order on $\La$ compatible with the partial order $\preceq$, the set
of objects $\O(\lambda)$
 indexed by $\Lambda$ with this order is
an  exceptional set generating $D^G(\Nt)$.

We now introduce another partial
 ordering $\leq$ on $\Lambda$. To this end, recall
the 2-sided Bruhat partial order on the 
affine Weyl group $W_{aff}$. For $\lambda\in \Lambda$ let
$w_\lambda$ be the minimal length representative of the coset
$W\lambda\subset W_{aff}$.
 We set
$\mu\leq \lambda$ if $w_\mu$ precedes $w_\lambda$ in the Bruhat
order.

We fix a complete order
%
 $\leq_{compl}$ on $\Lambda$ compatible
with $\leq$; we assume that $\{\mu \ |\ \mu \leq_{compl} \la \}$ is
finite for any $\la$.
We define the exceptional set $\nab_\lambda$ to be the $
\leq_{compl}$ mutation of the set $\O(\lambda)  $. It follows from
the above that $\nab^\la$ is an exceptional set generating
$D^G(\Nt)$. We define the {\em equivariant exotic $t$-structure} to
be the $t$-structure of the exceptional set $\nab^\la$, the objects
in the heart will be called equivariant exotic sheaves.

We now state compatibility between exotic and equivariant exotic
$t$-structures. Roughly speaking, over an orbit in $\N$ of
codimension $2d$ they differ by a shift by $d$. To state this
property more precisely,
 we need to recall the {\em perverse coherent
$t$-structure} \cite{izvrat}. Let $H$ be an algebraic group (assumed
for simplicity of statements connected) acting on an algebraic
variety $X$. Let $\bp$ be a function, called the perversity
function, from the set of $H$-invariant points of the scheme $X$
 to $\Zet$. We assume that $\bp$ is
{\em strictly 
 monotone and comonotone}, i.e. for points $x,\,y$, such that $x$ lies in the closure of $y$
we have $\bp(y) < \bp(x) < \bp(y)+\dim (y) -\dim (x)$. Then one can
define the perverse $t$-structure on $D^H(X)$, which
 shares some
properties with  perverse $t$-structure on the derived category of
constructible sheaves \cite{BBD}. For example, each perverse
coherent sheaf (i.e., object in the heart of the $t$-structure) has
finite length, and irreducible objects are in bijection with pairs
$(O,\LL)$, where $O\subset X$ is an $H$-orbit, and $\LL$ is an
irreducible $H$-equivariant vector bundle on $O$. In particular, if
the action is such that all orbits have even dimension, then the
perversity function $\bp(x)=\frac{\codim x}{2}$, called the middle
perversity, is strictly monotone and comonotone. It is well-known
that the adjoint action of a semi-simple group $G$ on the nil-cone
$\N$ has even dimensional orbits.

This construction works also for the category $D^H(A)$, where $A$ is
a coherent sheaf of associative $\O_X$ algebras equivariant under
$H$. 

\begin{proposition}\label{pervrA}
There exists a $G$-equivariant vector bundle $\E$ on $\Nt$, such
that $\E$, with the $G$-equivariant structure forgotten, is a
projective generator for the heart of the exotic $t$-structure.

We have an equivalence $\F\mapsto R\Hom(\E, \F)$ between $D^G(\Nt)$
and $D^G(A)$, where $A=End(\E)^{op}$. Under this equivalence the
equivariant exotic $t$-structure corresponds to the perverse
coherent $t$-structure of the middle perversity.
\end{proposition}


\subsection{Grading on exotic sheaves and canonical bases}

\subsubsection{Graded equivariant category
 and positivity by Frobenius weights}
We proceed to state a deep property of exotic sheaves related to an
additional grading on the Ext spaces between them. Recall the
current assumption that $char (\kk)=0$.

Consider the category $D^{G\times \Gm}(\Nt)$, where $\Gm$ acts on
$\Nt$ by $t:x\mapsto t^2 x$. For $d\in \Zet$ let $\F\mapsto \F(d)$
denote twisting by the $d$-th power of the tautological character of
$\Gm$.


We refer to \cite{Hum} for an elementary description of a canonical
lifting  $\tii \del_\la$, $\tii \nab^\la$ of $\del_\la$, $\nab^\la$
to $D^{G\times \Gm}$.
 This also fixes a lifting $\tii L$ of each irreducible
equivariant exotic sheaf $L$ to $D^{G\times \Gm}$.

\begin{theorem}\label{dpst}  For irreducible exotic equivariant
sheaves $L_1$, $L_2$ we have

 $\Ext^i(\tii L_1, \tii L_2(d))=0$ for
$d\leq 0$ and all $i$.

\end{theorem}


\begin{remark}
The Theorem follows from results of  \cite{B} on relation between
 exotic
sheaves and perverse sheaves on the affine flag manifold of the dual
group, see also Proposition \ref{qFr} below. They allow to deduce
the Theorem from Gabbers's Theorem \cite{BBD} on positivity of
weights of Frobenius action on Ext's between pure perverse sheaves
of the same weight. Thus it is the least elementary of the results
mentioned so far in this text.

The motivation for the Theorem is its consequence below, which shows
(in most cases) that classes of exotic sheaves form a {\em canonical
basis} in the Grothendieck group. This is parallel to the proof of
the Kazhdan-Lusztig conjecture: according to Soergel, cf.
\cite{SoerICM}, the latter is equivalent to the statement that for a
certain explicitly defined graded version of
Bernstein-Gel'fand-Gel'fand category $O$ the grading on $\Ext^1$
between irreducible objects has vanishing components of non-positive
degrees. The only known way to prove this vanishing is to identify
category $O$ with a category of perverse sheaves or Hodge
$D$-modules, and use Gabber's Theorem or its Hodge theoretic
analogue.
\end{remark}

\begin{remark}\label{cohotilt}
Another application of Theorem \ref{dpst} is explained in
\cite{Hum}. Together with {\em Koszul duality} formalism of
\cite{BGS} it allows to show that equivariant exotic sheaves control
cohomology of quantum groups at a root of unity with coefficients in
a tilting module.
\end{remark}




\subsubsection{Non-equivariant graded category and canonical
bases}\label{neqgr} We fix $X=\Nt$. Recall the category
$\A^f=\A_{A_0}^f\subset \A$ of exotic sheaves of finite length. It
is easy to see that $\A^f = \oplusl_{e\in \N} \A_e$, where
$\A_e=\A\cap \D_e$, and $\D_e\subset D(\Nt)$ is the full subcategory
of complexes whose cohomology sheaves are set-theoretically
supported on $\BB_e=\pi^{-1}(e)$. We have $K(\A_e)\cong K(\BB_e).$
Furthermore, the Chern character map provides an isomorphism
$K(\BB_e)_F\cong
H_*^{BM}(\BB_e)_F,$ 
 where $F$ denotes the coefficient
field of characteristic zero ($\Ce$ or $\Qlbar$), see, e.g,
 \cite{BMR}I.

The classes of irreducible objects form a basis in $K(\A_e)$. We
proceed to explain the properties of the category, which are needed
to relate this basis to {\em the canonical bases} in
$H_*^{BM}(\BB_e)$. The definition of the latter is due to Lusztig
\cite{LKthII}II, and follows the example of Kashiwara's
characterization of crystal bases \cite{Kashiwa}. More precisely,
Lusztig suggested a way to characterize a basis in
$H_*^{BM}(\BB_e)$, and conjectured that a basis satisfying his
axioms exists; he showed that it is then unique (up to a sign). We
will not recall Lusztig's characterization in detail; instead we
describe its structure and  explain the properties of exotic
sheaves, which imply (modulo a technicality, which is easy to check
in many cases)
 that Lusztig's axioms are satisfied by the basis
of irreducible exotic sheaves.

\medskip

One can find a homomorphism $\varphi:SL(2)\to G$, such that
$d\varphi$ sends the standard upper triangular generator of $sl(2)$
to $e$. Then we get an action $a_\varphi$ of the multiplicative
group $\Gm$ on $\gg$ given $a_\varphi(t):x\mapsto t^2\cdot ad
(\varphi (diag(t^{-1},t))) x$. This action fixes $e$.

 We let
$\D_e^{\Gm}\subset D^{\Gm}(\Nt)$ 
be the full subcategory of complexes, which are set theoretically
supported on $\pi^{-1}(e)$. Twisting by the tautological character
of $\Gm$ defines an auto-equivalence of this category, which we
denote by $\F\mapsto \F(1)$.
%
The exotic $t$-structure is inherited by the $\Gm$-equivariant
category; we let $\A_e^{gr}$ denote the heart of the latter. It is
easy to see that the forgetful functor $\D_e^{\Gm} \to \D_e$ sends
$Irr(\A_e^{gr})$ to  $Irr(\A_e)$, where $Irr$ stands for the set of
isomorphism classes of irreducible objects. This gives a bijection
$Irr(\A_e^{gr})/\Zet\cong Irr(\A_e)$, where $\Zet$ acts on
$Irr(\A_e^{gr})$ by $\F\mapsto \F(n)$. We also have
$K(\A_e^{gr})\cong K^{\Gm}(\BB_e)\cong K(\BB_e)[v,v^{-1}],$
where multiplication by $v$ corresponds to twisting by the
tautological character of $\Gm$.


 The canonical basis in $K(\BB_e)[v,v^{-1}]$ is characterized (up to a sign) by two
properties: {\em invariance under an involution} and {\em asymptotic
orthogonality} \cite{LKthII}II. These are reflected, respectively,
in categorical properties (i) and (ii) in the next Theorem.





Notice that the action of $B_{aff}$ on $\D_e$ is inherited by
$\D_e^{\Gm}$. Recall that $\SS$ is the Grothendieck-Serre duality.

In view of Theorem \ref{tstralc} and Proposition \ref{minusA}, the
contravariant auto-equivalence $\tii w_o \circ \SS$ is $t$-exact
with respect to the $t$-structure corresponding to the fundamental
alcove $A_0$, hence it permutes irreducible objects of $\A_{A_0}$;
here $w_o\in W$ is the long element.

\begin{theorem}\label{cordpst} 
There exists a canonical section of the map $Irr(\A_e^{gr})\to Irr
(\A_e)$, $L\mapsto \tii L$, such that

i) The image of the section is invariant  under every automorphism
of $G$ which is identity on the image of $\varphi _e$, and also
under $\tii w_o \circ \SS$.

ii) $\Ext^1_{\bar \A_e^{gr}}(\tii L_1, \tii L_2(i))=0$ for $i\leq 0$
and any $L_1, \, L_2\in Irr(\A_e)$; here $\bar \A_e^{gr}\subset
\A_e$ is the full subcategory of objects where the ideal of the
point $e$ in $\O(\gg)$ acts by zero.
\end{theorem}

{\em Comments on the proof.} The Theorem can be deduced formally
from Theorem \ref{dpst} and Proposition \ref{pervrA}. Thus its proof
relies on ideas of geometric Langlands duality used in \cite{AB},
and on Gabber's Theorem (see comments after Theorem \ref{dpst}).

\begin{corollary}\label{suppid}

Suppose that the involution $\tii \beta$ defined in \cite{LKthII}II,~\S 5.11
 induces
identity on the specialization at $q=1$. Then Conjecture 5.12
of {\em loc. cit.}, except, possibly, 5.12(g), holds;
moreover, the signed basis ${\bf B}_{\BB_e}^\pm$,
 whose existence is conjectured
in {\em loc. cit.},
 is formed by the classes of the objects $\tii
L$, where $L$ runs over  irreducible objects in $\A_e$.

\end{corollary}


\begin{remark}

The assumptions of the Corollary are easy to check in many cases,
e.g., if the nilpotent element $e$ is regular in a Levi subalgebra.
\end{remark}



\begin{remark} In fact, in \cite{LKthII}II Lusztig works with sheaves which are also equivariant under
a maximal torus in the centralizer of $e$. We omit this version here
to simplify notations, treating this set-up does not involve new
ideas.
\end{remark}

\begin{remark}
Validity of Conjecture 5.12(g) of \cite{LKthII}II is related to the
following
 question. Let $Y\subset \gg$ be a (Slodowy) transversal slice to a
$G$ orbit in $\N$, and $X=\Nt\times _\gg Y$. Let $A_X$ be as in
Theorem \ref{tstr}. One can show that $A$ can be endowed with a
natural grading; moreover, Theorem \ref{cordpst} is equivalent to
the fact that this grading can be chosen so that the graded
components of negative degree vanish, while the component of degree
zero is semi-simple. The question is whether the resulting graded
algebra is Koszul. If $e$ is subregular, then the positive answer is
easy to prove.
\end{remark}

\subsubsection{Independence of the (large) prime}
It is not hard to show that (co)homology of the Springer fiber is
independent of the ground field $\kk$, i.e. we have canonical
isomorphisms $H_\bu ^{BM}(\BB_e^\kk)\cong H_\bu^{BM}(\BB_e^\Ce)$,
where the upper index denotes the ground field, and $H_\bu^{BM}$
stands for $l$-adic Borel-Moore homology, $l\ne char(\kk)$.

 The definition of
the exotic $t$-structure is not specific to a particular ground
field.
 This allows one to prove the following.

\begin{proposition}\label{indep}
For all but finitely many prime numbers $p$ the following is true.
The classes in $H_\bu^{BM}(\BB_e^\kk)=H_\bu^{BM}(\BB_e^\Ce)$ of
irreducible exotic sheaves over $\kk$ of characteristic $p$ coincide
with the classes of irreducible exotic sheaves over $\Ce$.
\end{proposition}


\section{$D$-modules 
in positive characteristic and localization Theorem}\label{se3}
\subsection{Generalities on crystalline
 $D$-modules in positive characteristic}
\subsubsection{Definition and description of the center} Let
$X$ be a smooth variety over the field $\kk$.

The sheaf $\DD=\DD_X$ of {\em crystalline differential operators}
(or differential operators without divided powers, or PD
differential operators) on $X$ is defined as the enveloping of the
tangent Lie algebroid, i.e.,  for an affine open $U\subset X$ the
algebra $\DD(U)$   contains the subalgebra $\O$ of functions, has an
$\O$-submodule identified with the Lie algebra of vector fields
$Vect(U)$ on $U$, and these subspaces generate $D(U)$ subject to
relations $\xi_1\xi_2-\xi_2\xi_1=[\xi_1,\xi_2]\in Vect(U)$ for
$\xi_1, \xi_2\in Vect(U)$, and  $\xi \cdot f -f\cdot \xi =\xi(f)$
for $\xi\in Vect(U)$ and $f\in \O(U)$.

If $char(\kk)=0$, then $\DD_X$ is the  familiar sheaf of differential
operators. From now on assume that $\kk$ is of characteristic $p>0$.
Then $\DD_X$ shares some features with the characteristic zero case;
 for example,
$\DD_X$ carries an increasing filtration ``by order of a
differential operator'', and the associated graded $gr(\DD_X)\cong
\OO_{T^*X}$ canonically. On the other hand, some phenomena are
special to the characteristic $p$ setting. We have an action map
$\DD_X\to End(\OO_X)$,  which is not injective, unlike in the case
of a characteristic zero.
 For example, if $X={\mathbb A}^1=Spec(k[x])$, the
section $\de_x^p\ne 0$ of $\DD_X$ acts by zero on $\OO$. Also,
$\DD_X$ has a huge center; for example, if $X={\mathbb
A}^n=Spec(k[x_1,\dots,x_n])$, then $x_i^p$ and $\de_{x_i}^p$ are
readily seen to generate the center $Z(\DD_{{\mathbb A}^n})$ freely.
More generally, for any $X$ the center $Z(\DD_X)$
 is freely generated by elements
of the form $f^p$, $f\in \OO_X$ and $\xi^p-\xi^{[p]}$, $\xi \in
Vect_X$, where $\xi^{[p]}$ is the {\em restricted power} of the
vector field $\xi$; it is characterized by
$Lie_{\xi^{[p]}}(f)=Lie_\xi^p(f)$ for
 $f\in \OO_X$,
where $Lie$ stands for the Lie derivative. The center $Z(\DD_X)$ is
canonically isomorphic to the sheaf of rings $\OO_{T^*X\tw}$ where
the super-index $(1)$ stands for Frobenius twist.\footnote{Recall
that Frobenius twist of a variety $X$ over a perfect field $\kk$ is
defined to be isomorphic to $X$ as an abstract scheme, with the
$\kk$-linear structure twisted by Frobenius.
Not only $X\cong X\tw$ as abstract schemes, but also $X\cong X\tw$
as $\kk$-schemes, provided that $X$ is defined over $\bF_p$. For this
reason we will sometimes identify $X$ with $X\tw$ and omit Frobenius
twist from notation.}
%
 Thus $\DD_X$ can be considered as a
quasi-coherent sheaf of
 algebras on $T^*X\tw$.

\subsubsection{Azumaya property} Recall that an {\em Azumaya algebra
}
 on a scheme $X$
is a locally free sheaf $\AA$ of associative $\OO_X$ algebras, such
that the fiber of $\AA$ at every geometric point is isomorphic to a
matrix algebra. 
The following fundamental observation is due to Mirkovi\' c and
Rumynin, though a weak form of it can be traced to an earlier work
\cite{hren}.

\begin{theorem}\cite{BMR}I
 $\DD_X$ is an Azumaya algebra of rank $p^{2\dim (X)}$ on $T^*X^{(1)}$.
 \end{theorem}

See  \cite{OV}, \cite{BelKa} for generalizations and applications.


Recall  that two Azumaya algebras $\AA$, $\AA'$ are called equivalent (we then
write $\AA\sim \AA'$) if they are Morita equivalent, i.e. if there
exists a coherent locally projective sheaf $\MM$ of $\AA-\AA'$
bimodules, such that
 $\AA'\iso End(\MM)^{op}$; we will then say that $\MM$ provides
an equivalence between $\AA$ and $\AA'$.
 In particular, an Azumaya algebra $\AA$ is {\em split} if $\AA\sim \OO_X$;
this happens iff $\AA\cong  End(\EE)$ for a vector bundle $\EE$.
 For
two equivalent Azumaya algebras $\AA$, $\AA'$ we have an equivalence
of categories of modules $Coh(X,\AA)\cong Coh(X,\AA')$, depending on
the choice of a bimodule providing the equivalence between $\AA$ and
$\AA'$; in particular, for a split Azumaya algebra we have
$Coh(X,\AA)\cong Coh(X)$.


For a smooth variety $X$ over a positive characteristic field, the
Azumaya algebra $\DD_X$ is not split unless $\dim (X)=0$. However,
it is split on the zero section, see \cite{OV} for more information.

We will also need a twisted  version of differential operators. If
$\LL$ is a line bundle on $X$, then one can consider the sheaf
$\DD^\LL=\DD_X^\LL$ of differential operators in $\LL$.
 A similar argument shows
that this is also an Azumaya algebra over $T^*X\tw$; moreover, we
have a canonical equivalence
\begin{equation}\label{twisteq}
\DD_X\sim D^\LL_X, \end{equation}
 given by the bimodule
$D_X\otimes_{\O(X)} \LL^{-1}$.

\begin{remark}\label{praz}
 Notice that if $\LL=\LL_0^p=Fr^*(\LL_0)$ for some line
bundle $\LL_0$, then we have a canonical isomorphism $\DD^\LL_X\cong
\DD_X$; however, the above equivalence $\DD^\LL_X\cong \DD_X$ is not
identity, but rather tensor product
 {\em over the
ring $\O_X^p=\O_{X\tw}$} with the line bundle $\LL_0\tw$.
\end{remark}

\subsection{Crystalline operators on $\BB$}
We now consider $X=\BB$. We abbreviate $\DD_\BB^{\O(\la)}=\DD^\la$.

\subsubsection{Splitting the Azumaya algebra}\label{split_sect}
It was mentioned above that  $\DD_\BB$ splits on the zero section.
In fact, we have the following stronger statement.

\begin{theorem}\label{spl}


a)  There exists an Azumaya algebra $\AA$ on $\N\tw$, such that
$\DD^{-\rho}\cong \pi^*(\AA)$.

b) For any $\la$ we have an equivalence of Azumaya algebras on $\N\tw$:
$\DD^\la\sim \pi^*(\AA)$.

c) $\DD^\la(\BB)$ is split on the formal neighborhood of every fiber of
$\pi$.
\end{theorem}
{\em Sketch of proof.}  (a) reduces to
irreducibility of baby Verma modules with highest weight $-\rho$, which
follows from \cite{BrGo}.
 It implies, moreover,
that the statement holds for $\A$ being the quotient of the
enveloping algebra  $U(\gg)$ by the central ideal corresponding to
$-\rho$.
 (b) follows from (a) in view of the equivalence
\eqref{twisteq}. Finally, (c) follows from (b), since every Azumaya
algebra over a complete local ring with an algebraically closed
residue field  splits.

\medskip

Let $\DD^\la-\rmod^f\subset Coh(\Nt\tw, \DD^\la)$
  the full subcategory
of sheaves, whose support (which is a subvariety in $\Nt\tw$) is
proper. Let $Coh^f(\Nt)\subset Coh(\Nt)$ be the full subcategory of
sheaves with proper support.

\begin{corollary}\label{eqccoh}
For every $\la\in \La$ we have an equivalence $ \DD^\la-\rmod^f
\cong Coh^f(\Nt)$.
\end{corollary}

\proof Since the target of $\pi$ is affine, a subscheme $Z$ in
$\Nt\tw$ is proper iff it lies in a finite union of nilpotent
neighborhoods of Springer fibers. Thus the claim follows from
Theorem \ref{spl}(b).

\medskip

For each $\la\in \La$ and $e\in \N$
we fix  the splitting bundle $\EE^\la_e$ for $\DD^\la$ on the formal
neighborhood of $\pi^{-1}(e)$ as follows. For $\la=-\rho$ we let $\EE_e^\la$
be the pull-back under $\pi$
of a splitting bundle for the Azumaya algebra $\A$ on the formal neighborhood of $e$ in
$\N\tw$. For a general $\la$ we get $\EE^\la$
from $\EE^{-\rho}$ by applying the canonical equivalence \eqref{twisteq}
between $\DD^{-\rho}$ and $\DD^\la$; thus $\EE^\la=\E^{-\rho}\otimes _{\O_\BB}
\O(\la+\rho)$.

We let $\Feq_\la$ denote the resulting equivalence between
$\DD^\la-\rmod^f$ and $Coh^f(\Nt)$. Notice that for $\la'=\la+p\mu$
the sheaves of algebras $\DD^\la$ and $ \DD^{\la'}$ are canonically
identified; however, the equivalences $\Feq_\la$ and $\Feq_{\la'}$
are different, cf. Remark \ref{praz}.



\subsubsection{Derived localization in positive
characteristic}\label{derloc} Let $U=U(\gg)$ be the enveloping
algebra.

Assume first that $char(\kk)=0$. Recall the famous Localization
Theorem \cite{BeBe}, \cite{BrylKa}, which provides an equivalence
$U^\la-\rmod\cong \DD^\la-\rmod(\BB)$, where $\la$ is a dominant
integral weight, $\DD^\la-\rmod$ denotes the corresponding twisted
$D$-modules category, and $U^\la-\rmod$ is the category of
$\gg$-modules with central character corresponding to $\la$. For two
 integral weights $\la, \, \mu$ the categories $\DD^\mu-\rmod$ and
$\DD^\la-\rmod$ can be identified by means of the equivalence
$T^\la_\mu:\F \mapsto \F\otimes \O(\la-\mu)$. If $\la, \, \mu$ are
dominant, then the global sections functors intertwine this
equivalence with  the translation functor, which provides an
equivalence $U^\la-\rmod\cong U^\mu-\rmod$.

Assume now that $\mu$ is integral regular, thus
$\mu=w(\la+\rho)-\rho$ for some dominant integral $\la$, $w\in W$.
Then the functor of global sections on $D_\mu-\rmod$ is no longer
exact; however, it follows from \cite{BeBeder} that the derived
functor $R\Ga=R\Ga_\mu:D^b(\DD^\mu-\rmod)\to D^b(U^\la-\rmod)$ is
still an equivalence. The triangle formed by the three equivalences
$R\Ga_\mu$, $T^\la_\mu$, $R\Gamma_\la$ does not commute. Thus we get
an auto-equivalence  $R_w$ of $D^b(\DD^\la-\rmod)$,
$R_w=R\Ga_\la^{-1} \circ R\Ga_\mu\circ T_\la^\mu$. In \cite{BeBeder}
it is shown that $R_w$ can be described by an explicit
correspondence, which makes it natural to call $R_w$ the {\em Radon
transform}, or {\em the intertwining functor}. Moreover, the
assignment $\tii w\mapsto R_w$ extends to an action of the Artin
braid group $B$ attached to $G$ on $D^b(U^\la-\rmod)$.

A part of this picture can be generalized to characteristic $p$.

The obvious characteristic $p$
 analogue of the above equivalence of abelian categories
 does not hold for any integral
$\la$. Indeed, it is well known that for any coherent sheaf $\F$ on
the (Frobenius twist  of) a smooth variety the sheaf $Fr^*(\F)$
carries a flat connection; in particular, so does the sheaf
$Fr^*(L)=L^{\otimes p}$, where
 $L$ is a line bundle.
 Thus for $\F\in \DD^\la-\rmod$ we have $\F\otimes L^p\in
\DD^\la-\rmod$. If $L$ is anti-ample and the support of $\F$ is
projective of positive dimension, then  some of the higher derived
functors $R^i\Gamma(\F\otimes L^{dp})\ne 0$ for large $d$.

However, we do have an analogue of the "derived" localization
Theorem. From now on assume that $char(\kk)=p>0.$

The center $Z$ of $U$ contains the subalgebra $Z_{HC} =U^G\cong
Sym(\hh)^W$, which we call the Harish-Chandra center. We have a
natural map $\La/p\La \to \hh^*/W$, $\la \mapsto d \la\mod W$. Thus
every $\la\in \La$ defines a maximal ideal of $Z_{HC}$. We let
$U^\la=U\otimes_{Z_{HC}} \kk$ denote the corresponding central
reduction. Notice that the set of weights  $\mu\in\La$, such that
the quotients $U^\la$ and $U^\mu$ of $U$ coincide,
 is precisely the $W_{aff}$-orbit of $\la$ with respect
to the  action $w\bu \la = p\, w(\frac{\la+\rho}{p})-\rho$.
We will say that $\la\in \La$ is $p$-regular if the stabilizer in $W$
of $\la+p\La\in \La/p\La$ is trivial.

We also have another central subalgebra $Z_{Fr}\subset U$, called
the Frobenius center. It is generated by expressions of the form
$x^p-x^{[p]}$, $x\in \gg$, where the {\em restricted power} map
$x\mapsto x^{[p]}$ is characterized by $ad(x^{[p]})=ad(x)^p$. Thus
maximal ideals of $Z_{Fr}$ are in bijection with points of
$\gg^*\cong \gg$.

Let $U^\la-\rmod$ denote the category of finitely generated $U^\la$-modules,
and let $U^\la-\rmod^f\subset U^\la-\rmod$ be the full subcategory of finite
length modules.

For a pair $\la\in \La$, $e\in \gg^*$  let $U^\la_{\hat e}-\rmod$ be
the category of finitely generated $U^\la$-modules, which are killed
by some power of the maximal ideal of $e$ in $Z_{Fr}$. This category
is zero unless $e\in \N$. We also have $U^\la-\rmod^f =\oplusl_{e\in
\N} U^\la_{\hat e}-\rmod$.

Let $\DD^\la_{\hat e}-\rmod \subset \DD^\la-\rmod$ be the full
subcategory of objects which are supported on a nilpotent
neighborhood of $\pi^{-1}(e)$; here we think of $\DD^\la$ modules as
sheaves on $\Nt\tw$ with an additional structure.

\begin{theorem}\label{BMRth} \cite{BMR}I a) For every $\la\in \La$
we have a natural isomorphism $\Ga(\DD^\la)\cong U^\la$.

b) If $\la\in \La$ is $p$-regular, then the derived global sections
functor provides an equivalence $R\Ga_\la: D(\DD^\la)  \iso
D(U^\la)$. It restricts to equivalences $D^b(\DD^\la\rmod^f)\cong
D^b(U^\la-\rmod^f)$, $D^b(\DD^\la_{\hat e}-\rmod)\cong
D^b(U^\la_{\hat e}-\rmod)$.
\end{theorem}

\begin{remark}\label{dgRem}
This Theorem has several versions and generalization. One can work
with the more general categories of twisted $D$-modules, thereby
obtaining a category of modules over an Azumaya algebra on the
formal
 neighborhood of $\Nt$
in $\gt$, or more general subschemes or formal completions of $\gt$.
For singular weights $\la$ there is a version of the Theorem that
relates derived categories of modules to sheaves on (the
neighborhoods of)  {\em parabolic Springer fiber} \cite{BMR}II.
Another construction works with differential operators on a partial
flag variety $G/P$ for a parabolic subgroup $P\subset G$, {\em loc.
cit.}, cf. also subsection \ref{parvar} above.

For a scheme $Y$ mapping to $\gg$ and satisfying the Tor vanishing
conditions of Theorem \ref{bract} we have an equivalence between the
derived category of modules over  Azumaya algebras on $\tii \bY$,
$\tii Y$ obtained as pull-back of the algebra of (twisted)
differential operators and derived category of modules over the
algebra of global sections.

If $Y$ is a transversal slice to a nilpotent orbit, then the algebra
of global section of the Azumaya algebra on $\tii Y$ is probably
related to Premet's quantization of Slodowy slices and generalized
Whittaker $D$-modules, see \cite{Premet}, \cite{GinGan}.

There exists a generalization of this result for $Y$ not satisfying
the $\Tor$ vanishing condition. It involves coherent sheaves on the
{\em differential graded} scheme, which is the derived fiber product
of $Y$ and $\gt$ over $\gg$. The particular case $Y=\{0\}$ is
closely related to the description of the derived category of the
principal block in representations of a quantum group at a root of
unity provided by \cite{ABG}.

\end{remark}


\begin{proposition}
a) A weight $\la\in \La$ is $p$-regular iff $\frac{\la+\rho}{p}$ lies
in some alcove.

b) The $t$-structure on $\DD^0-\rmod$ induced by the equivalence
$R\Ga^\la\circ T_0^\la$ for a $p$-regular $\la$ depends only on the
alcove of $\frac{\rho+\la}{p}$ (see beginning of section
\ref{derloc} for notation).
\end{proposition}

Thus we get a collection of $t$-structure on $\DD^0-\rmod$ indexed
by
 alcoves; we denote the $t$-structure attached to $A\in \Alc$ by $D^{<0}_A(\DD)$, $D^{\geq 0}_A(\DD)$.
The following properties of the collection follow from \cite{BMR}.


\begin{theorem}\label{abcd}
a) Let $A_1$, $A_2$ be two alcoves. If $A_1$ lies above $A_2$, then
$D^{> 0}_{A_1}(\DD)\supset D^{> 0}_{A_2}(\DD)$.

b) There exists an action of $B_{aff}$ on $D(\DD)$, such that the following
holds.
Let $\la\in \La$, $\frac{\la+\rho}{p}=w(\frac{\rho}{p})$ for $w\in W_{aff}$,
and let $A$ be the alcove of $\frac{\la+\rho}{p}$.
Then $R\Ga_\la\cong R\Ga_0\circ b_{A_0,A_1}$.
Thus $b_{A_0,A_1}$ sends the $t$-structure
 $D^{<0}_{A_0}(\DD)$, $D^{\geq 0}_{A_0}(\DD)$ to the $t$-structure
 $D^{<0}_A(\DD)$, $D^{\geq 0}_A(\DD)$.

c) The restriction of the $B_{aff}$ action to $D^b(\DD^\la-mod^f)
\cong D^b(Coh^f(\Nt))$ coincides with the restriction of the
action from Theorem \ref{bract}.

d) The derived global sections functor $R\Gamma$ on
$D^b(Coh^f(\Nt))$ is $t$-exact with respect to the $t$-structure
induced from $( D^{<0}_{A_0}(\DD)$, $D^{\geq 0}_{A_0}(\DD))$ via the
equivalence $\Feq_0$ (see the end of section \ref{split_sect} for
notation).

\end{theorem}

\begin{corollary}\label{AU}
The $t$-structure on $Coh^f(\Nt)$ induced by the exotic
$t$-structure coincides with the one induced from
$(D^{<0}_{A_0}(\DD)$, $D^{\geq 0}_{A_0}(\DD))$ via the equivalence
$\Feq_0$ (see the end of \ref{split_sect}).

In particular, we have a Morita equivalence
 $\A_e\sim U-\rmod_{\hat
e}^\la$ for every $p$-regular $\la\in \La$.
\end{corollary}

The Corollary follows by comparing Theorem \ref{abcd} with Theorem
\ref{tstralc}.

\begin{corollary}\label{Ue}
a) Let $U_e^\la$ denote the specialization of the enveloping algebra
$U(\gg)$ at the central character corresponding to  $e\in \N$ and a
regular integral weight $\la$. Then we have a canonical isomorphism
$K(U_e^\la)_F\cong H_\bu^{BM}(\BB_e)_F$, where $F$ is a field of
characteristic zero.

b) The images of the irreducible modules under this isomorphism is
independent of the base field $\kk$, except for a finite number of
values of the characteristic.

\end{corollary}

Part (a) of the Corollary follows directly from Theorem \ref{BMRth}
(cf. the discussion preceding Proposition \ref{pospos}), while part
(b) follows from  Corollary \ref{AU} and Proposition \ref{indep}.

\begin{remark}
For $e=0$ part (a) of the  Proposition is standard, and part (b) can
be deduced from \cite{AJS}. Our method uses the principal tool of
\cite{AJS}, namely, the reflection functors, in the disguise of the
braid group action; geometry of the Springer map is the new
ingredient.
\end{remark}

In fact, we have the following stronger, though more difficult
statement. We will say that a basis in $H_\bu^{BM}(\BB_e)$  is
canonical if it is the image of a basis in the equivariant
Grothendieck group $K^{\Gm}(\BB_e)$, satisfying Lusztig's axioms
\cite{LKthII}II,  under forgetting the equivariance composed with
the Chern character map. According to a result of \cite{LKthII}II
such a basis is unique up to multiplication of some of its elements
by $-1$, if it exists.

\begin{corollary}  Enforce the assumption of Corollary \ref{suppid}.
Then for almost all $p=char(\kk)$ the isomorphism
$K^0(U-mod_e^0)_F\cong H_*^{BM}(\BB_e)_F$ of Corollary \ref{Ue}(a)
sends classes of irreducible objects to elements of a canonical
basis. Thus Conjecture 17.2 of \cite{LKthII}II holds in this case.
\end{corollary}

This Corollary is immediate from Corollary \ref{suppid}  together
with Theorem \ref{abcd}(d). Thus its proof, unlike the proof of
Corollary \ref{Ue}, relies on Gabber's Theorem \cite{BBD} and ideas
of local geometric Langlands duality, on which the results of
\cite{AB} are based.

\begin{remark}
The particular case $e=0$ of the Conjecture 17.2 of \cite{LKthII}II
is well known to imply the previous Lusztig conjectures
\cite{LConj}, which describe characters of algebraic groups in
finite characteristics. Lusztig's program for a proof of these
conjectures has been carried out by several authors. An alternative
proof is given in \cite{ABG}.

Notice that the strategy of proof of this conjecture outlined above
does not use quantum groups.
\end{remark}





\section{Perverse sheaves on  affine flags of the dual group
(local geometric Langlands).}\label{se5}


\subsection{Generalities on geometric Langlands duality}
 Recall that $\LG$ is the
group dual to $G$ in the sense of Langlands. Several good
 surveys of geometric Langlands duality program has appeared
recently \cite{Edik1}, \cite{Edik}, \cite{GaiICM},
 so I will only
briefly recall the set-up.

The geometric Langlands duality is a categorification of the classical
Langlands duality for functional fields. The latter seeks to attach
an automorphic form to a homomorphism from a version of the
 Galois group to $\LG$. In other words, the problem
is to provide
a spectral decomposition for Hecke operators acting in the space of
 automorphic functions, and relate the space of spectral parameters
to homomorphisms of the Galois group  to the dual group. As was
probably first observed by A.~Weil, 
 in the case of
 a functional field
the automorphic space in question is the set of isomorphism classes
of $G$-bundles, possibly with an additional level structure, on an
algebraic curve over $\Fq$.  Thus it is the set of $\Fq$ points of
the corresponding moduli space (stack).

Passage to the geometric duality theory is based on the following
variation of Grothendieck's sheaf-function correspondence principle.
The variation says that for an algebraic variety (or stack) over
$\Fq$ a natural categorification of the space of functions on the
set $X(\Fq)$ is the derived category of $l$-adic sheaves on $X$.
Thus the objective of the geometric duality theory is a spectral
decomposition of the  derived category of $l$-adic sheaves on a
moduli space of $G$-bundles, where the space of spectral parameters
is identified with the space of $\LG$ local systems. It is a
non-trivial, and not completely solved, problem to assign a formal
meaning to the previous sentence; however, in some cases it amounts
to an equivalence between the $l$-adic derived category of the
moduli stack and the derived category of coherent sheaves on a stack
mapping to the stack of local systems.

The above formulations referred to a more developed {\em global}
version of the theory. However, the classical Langlands conjectures
have both a global and a local version. The global one provides a
conjectural classification of automorphic representations of the
group of adele points of a reductive group over a global field, i.e.
either a number field, or the field of rational functions on a curve
over a finite field. The local one describes all irreducible
representations of a reductive group over a local field; recall that
a functional local field is a field of formal Laurent series
$\Fq((t))$. The geometric theory studies the derived category of
$l$-adic sheaves on
 homogeneous spaces of the {\em formal loop group} $\bLG((t))$ of the
 dual group $\LG$. The latter is a group ind-scheme over $\Fq$,
 whose group of $\Fq$ points is identified with $\LG(\Fq((t)))$.
The results are expected to link such $l$-adic derived categories to
coherent sheaves on spaces related to $G$-local systems on the
puncture formal disc, cf. \cite{Edik_Denis}.
%

\subsection{Results of \cite{AB}, \cite{B}, \cite{Betobe}}
Some particular results of this type have been achieved in {\em loc.
cit.} 

\subsubsection{Statement of a result}
Recall that the Iwahori subgroup $I\subset \LG((\Fq(t)))$ consists
of those maps from the punctured formal disc to $\LG$, which can be
extended to a map from the whole disc, so that the image of the
closed point is contained in a fixed Borel subgroup $\LB\subset
\LG$.

We have a  group subscheme (pro-algebraic group) $\bI\subset
\bLG((t))$, such that  $\bI(\Fq)=I$. The affine flag space $\Fl$ of
$\LG$ is the homogeneous space $\bLG((t))/\bI$. It is  an
ind-algebraic variety such that $\Fl(\Fq)=\LG(\Fq((t)))/I$.
 The
group $\bI$ acts on $\Fl$. The orbits of this action, called affine
Schubert cells, are in bijection with the affine Weyl group
$W_{aff}$.

Let $\P$ denote the category of perverse sheaves on $\Fl$, which are
equivariant with respect to the prounipotent radical of $\bI$. Let
$\P^I\subset \P$ be the full subcategory of $\bI$ equivariant
sheaves.

Let $^{nf}\P^I\subset \P^I$ be the Serre subcategory generated by
irreducible objects, corresponding to those $w\in W_{aff}$, which
are not the minimal length representatives of a left $W$ coset. Let
$^f\P=\P/^{nf}\P$ be the Serre quotient category.

\begin{remark}
To clarify the definition of $^f\P$ we remark  that this category
can be also described as the category of Iwahori-Whittaker sheaves
\cite{AB}. Thus it is related to the Whittaker model, which is one
of the main tools in representation theory of reductive groups over
local and global fields.
\end{remark}

\begin{theorem}\label{equiv}
a) \cite{Betobe} We have a canonical equivalence $D^b(\P)\cong
D^G(\gt\times _\gg \Nt)$.

b) \cite{AB} We have a canonical equivalence $D^b(^f\P)\cong
D^G(\Nt)$. The image of $^f\P$ under this equivalence consists of
equivariant exotic sheaves.

\end{theorem}

The Theorem is motivated by the  known isomorphisms of
Grothendieck groups; the question of possibility of such (or
similar) equivalence has been raised, e.g., by V.~Ginzburg, see
Introduction to \cite{CG}. More precisely, the Grothendieck group of
the two categories appearing in Theorem \ref{equiv}(a) are
isomorphic to the
 group algebra of the affine Weyl group of $\LG$. A more
 interesting version of the isomorphism is obtained by replacing the
 categories by their graded version: $D^G(\gt\times _\gg \Nt)$ is
 replaced by $D^{G\times \Gm}(\gt\times \Nt)$, while the definition
 of the graded version of $\P$ is more subtle
 (cf. Proposition \ref{qFr} below
 and also \cite{BGS}). The corresponding Grothendieck groups turn
 out to be isomorphic to {\em the  affine Hecke algebra},
 see \cite{CG}, \cite{LKthII}I.

Similarly, the  Grothendieck groups of both categories appearing in
Theorem \ref{equiv}(b) are identified with the {\em anti-spherical}
module over the extended affine Weyl group, while in the graded
version of the theory we get the anti-spherical module over the
affine Hecke algebra. Here by the anti-spherical module we mean the
induction of the sign representation from the finite Weyl group
(respectively, Hecke algebra) to the affine one.

I would like to emphasize that this isomorphism of the two
realizations of the affine Hecke algebra is the key step in the
proof of classification of its representations due to Kazhdan and
Lusztig \cite{KL} (see also \cite{CG}), which establishes a
particular case of the local Langlands conjecture. This is another
illustration of the relation of Theorem \ref{equiv} to local
Langlands duality.

\medskip

The proof of Theorem \ref{equiv} builds on previously known
constructions of categories related to $G$ in terms of perverse
sheaves on homogeneous spaces for $\bLG((t))$. The first important
result is {\em the geometric Satake isomorphism} \cite{Ginzb},
\cite{MV}, \cite{BD}, which identifies the tensor category $Rep(G)$
of algebraic representations with the category of perverse sheaves
on the affine Grassmannian $\Gr=\bLG((t))/\GO$ equivariant with
respect to $\GO$. Here $\GO\subset \bLG((t))$ is the group
subscheme, such that $\GO(\Fq)$ consists of maps which extend to the
non-punctured disc. Furthermore, Gaitsgory \cite{Zen} used this
result to provide a categorification of the description of the {\em
center} of the affine Hecke algebra. Using some ideas of
I.~Mirkovi\' c we observe that the so-called Wakimoto sheaves
provide a categorification of the maximal abelian subalgebra in the
affine Hecke algebra due to Bernstein, see, e.g., \cite{CG},
\cite{LKthII}I.
 The maximal
projective object in the category of sheaves on the finite
dimensional flag variety  of $\LG$ smooth along the Schubert
stratification, which plays a central role in Soergel's description
of category $O$, cf. \cite{SoerICM}, is
  a categorification of the
$q$ anti-symmetrizer (an element of the finite Hecke algebra, which
acts by zero in all irreducible representation except for the sign
representation). Under the equivalence of Theorem \ref{equiv}(a) it
corresponds to the structure sheaf of $\gt\times _\gg \Nt$. A
combination of these ingredients yields a proof of the Theorem.

\subsubsection{Possible generalizations}
It is natural to ask if the multiplication in the affine Hecke
algebra corresponds to a monoidal structure on the derived
categories of coherent sheaves and constructible sheaves appearing
in Theorem \ref{equiv}(a). In order to get such a monoidal
structure, we need to replace the categories defined above  by
closely related ones with the same Grothendieck group. One way to do
it is as follows. Let $\bI'$ be the pro-unipotent radical of $\bI$.
Let $\P'$ be the category of perverse sheaves on "the basic affine
space" $\bLG((t))/\bI'$, which are $\bI$-monodromic with unipotent
monodromy. Then convolution provides the derived category $D^b(\P')$
with a monoidal structure. Notice that this monoidal category does
not have a unit object, though this can be repaired by adding some
pro-objects to the category, the unit object is then the free
pro-unipotent local system on $\bI/\bI'\subset \bLG((t))/\bI'$.

On the dual side we consider the category $D^G(\gt\times _\gg\gt)$.
Then one can show that convolution provides this category with a
monoidal structure. Let $Coh^G(\gt\times _\gg\gt)'\subset
Coh^G(\gt\times_\gg \gt)$ denote the full subcategory of complexes,
whose cohomology sheaves are set-theoretically supported on the
preimage of $\N\subset \gg$. A standard argument shows that it
yields a full embedding of derived categories $D^G(\gt\times
_\gg\gt)':=D^b (Coh^G(\gt\times _\gg\gt)')$ into  $D^G(\gt\times_\gg
\gt)$. The full subcategory $D^G(\gt\times _\gg \gt)'\subset D^G(\gt
\times _\gg \gt)$ is closed under the convolution product, though it
does not contain the unit object $\delta_*(\O)$, where
$\delta:\gt\to \gt\times _\gg \gt$ is the diagonal embedding.

It is easy to see that the push-forward (respectively, pull-back)
functors $Coh^G(\gt\times _\gg \Nt)\to Coh^G(\gt\times_\gg \gt)'$,
$\P\to \P'$ induce isomorphisms of Grothendieck groups.

\begin{theorem} \cite{Betobe}
We have a natural monoidal equivalence $D(\P')\cong D^G(\gt\times
_\gg \gt)'$.
\end{theorem}


\begin{remark}
Another version of Theorem \ref{equiv} links the monoidal $\bI$
equivariant derived category to the monoidal derived category of
$G$-equivariant coherent sheaves on the fiber square of $\Nt$ over
$\gg$. An additional subtlety in this case is related to
nonvanishing of  $\Tor_{>0}^{\O_\gg}(\O_\Nt, \O_\Nt)$. One actually
has to take these $\Tor$ groups into account by working with the
{\em derived fiber product}, which is a differential-graded scheme,
rather than an ordinary scheme. This issue does not arise in the
other settings mentioned above, because
$\Tor_{>0}^{\O_\gg}(\O_\gt,\O_\gt)=0=\Tor_{>0}^{\O_\gg}(\O_\gt,\O_\Nt)$.
However, one has to work with differential graded schemes in order
to define the convolution product on $D^G(\gt\times _\gg\gt)$.
\end{remark}

\subsubsection{Relation to the material of section \ref{se2}}
Many of the constructions from section \ref{se2} are motivated by
the equivalences of Theorem \ref{equiv}.

For example, the categories $D(\P)$, $D(\P')$ carry a natural
$B_{aff}$ action by {\em Radon transforms}, cf. beginning of section
\ref{derloc}, where a similar structure for a finite dimensional
flag variety is mentioned. To define the action we recall that the
$\bLG((t))$ orbits on $\Fl^2$ are indexed by the affine Weyl group.
If $\Fl_w^2$ is the orbit corresponding to $w\in W_{aff}$, and
$pr_i^w:\Fl_w^2\to \Fl$ are the projections, where $i=1,2$ then we
define a functor $R_w:D(\P)\to D(\P)$ by $R_w(\F)=pr^w_{2*}
pr_1^{w*}(\F)$. Then we have an action of $B_{aff}$ on $D(\P)$,
$D(^f\P)$, such that $\tii w\mapsto R_w$. Under the equivalences of
Theorem \ref{equiv}(b) this action corresponds to the action
described in section \ref{se2}.


Finally, I would like to quote the statement that allows to link the
grading on $\Ext$ spaces appearing in  Theorem \ref{dpst} to
Frobenius weights, thus providing a way to prove Theorem \ref{dpst}.
To state it we introduce the following notation.
Let $\Phi$ be either of the two equivalences appearing in Theorem
\ref{equiv}. Let $Fr$ be the autoequivalence of the corresponding
category of constructible sheaf, sending a sheaf to its pull-back
under the Frobenius morphism.
 Let $\bq$ be an automorphism of either $\Nt$ or $\gt\times _\gg \Nt$
given by $(\b,x)\mapsto (\b, qx)$ or $(\b_1,\b_2,x)\mapsto (\b_1,
\b_2, qx)$ respectively; here $q$ stands for the cardinality of the
base finite field $\Fq$.

\begin{proposition}\label{qFr} (cf. \cite{AB})
 We have a canonical isomorphism $\Phi\circ
\bq^*\cong Fr\circ \Phi$.
\end{proposition}

\end{document}